\documentclass[12pt,a4 paper]{amsart}

\pdfoutput=1

\usepackage{t1enc}
\usepackage[english]{babel}
\usepackage{amsmath}
\usepackage{amsfonts}
\usepackage{amssymb}
\usepackage{amsthm}
\usepackage{graphicx}
\usepackage{xcolor}
\usepackage{lmodern}
\usepackage[all]{xy}

\newtheorem{thm}{Theorem}[section]
\newtheorem{cor}{Corollary}[section]
\newtheorem{prop}{Proposition}[section]

\newtheorem*{nndef}{Definition} 
\newtheorem*{nnthm}{Theorem} 
\newtheorem*{nnexams}{Examples} 
\newtheorem*{nnrem}{Remark} 
\newtheorem*{nnsrem}{Side remarks} 

\DeclareMathOperator{\C}{\mathbb C}
\DeclareMathOperator{\R}{\mathbb R}
\DeclareMathOperator{\Z}{\mathbb Z}
\DeclareMathOperator{\N}{\mathbb N}
\DeclareMathOperator{\E}{\mathbb E} 

\DeclareMathOperator{\quat}{\mathbb H} 

\DeclareMathOperator{\qps}{{\rm P}^{n-1}(\mathbb H)} 

\DeclareMathOperator{\m}{\rm M} 
\DeclareMathOperator{\glg}{\rm GL} 
\DeclareMathOperator{\slg}{\rm SL} 
\DeclareMathOperator{\sog}{\rm SO} 
\DeclareMathOperator{\org}{\rm O} 
\DeclareMathOperator{\sug}{\rm SU} 
\DeclareMathOperator{\ung}{\rm U} 
\DeclareMathOperator{\spg}{\rm Sp} 
\DeclareMathOperator{\sla}{\mathfrak{sl}} 
\DeclareMathOperator{\spa}{\mathfrak{sp}} 

\DeclareMathOperator{\harm}{\rm H} 
\DeclareMathOperator{\sph}{\mathcal{Y}} 

\DeclareMathOperator{\tra}{{\rm{Tr}}} 

\DeclareMathOperator{\pild}{\pi_{ \textit{i} \lambda,\delta}} 

\newcommand{\trans}{{}^t\!}

\begin{document}

\title[Special functions, K-types and symplectic groups]{Special functions associated\\
with $K$-types of degenerate\\
principal series of $\spg(n,\C)$}

\author{Gr{\'e}gory Mendousse}

\date{}

\maketitle

\begin{abstract}
\
\noindent We study $K$-types of degenerate principal series of $\spg(n,\C)$ by using two realisations of these infinite-dimensional representations. The first model we use is the classical compact picture; the second model is conjugate to the non-compact picture via an appropriate partial Fourier transform. In the first case we find a family of $K$-finite vectors that can be expressed as solutions of specific hypergeometric differential equations; the second case leads to a family of $K$-finite vectors whose expressions involve Bessel functions.



\end{abstract}

\section{Introduction} \label{section-introduction}
\noindent Special functions play a prominent role in mathematics and physics. They are solutions of important equations and come in various forms, in terms of series or integrals. We are interested in the way special functions connect to representation theory. That such a connection exists is well known. Special functions appear in the study of spherical functions (see e.g. \cite{Helgason1984}, Chapter IV), of matrix coefficients of representations (see e.g. \cite{Vilenkin1993-1}, Section 4.1) and of symmetry breaking operators (see \cite{Kobayashi-Kubo-Pevzner2016} and \cite{Kobayashi-Pevzner2016}). In this work, we focus on the $K$-types of  degenerate principle series of $\spg(n,\C)$, looking for $K$-finite vectors that have explicit formulas in terms of special functions. Here, $K$ refers to $\spg(n)$.
\vskip 8pt
\noindent The reason why we choose to work with degenerate principle series of $\spg(n,\C)$ is twofold.
\vskip 8pt
\noindent First, the geometric setting we choose is motivated by the following sequence of groups of isometries of vector spaces over different number fields:
$$\xymatrix{
\spg(n) \ar@{~>}[d] & \subset & \sug(2n) \ar@{~>}[d] & \subset & \sog(4n) \ar@{~>}[d] \\
\quat^n & \simeq & \C^{2n} & \simeq & \R^{4n}.}$$
This sequence enables us to use a refined version of the classical theory of spherical harmonics. Moreover, the non-commutativity of the skew field $\quat$ of quaternions adds a rich ingredient to the underlying representation theory and harmonic analysis: see for example \cite{Howe-Tan1993},  \cite{Kobayashi1992},  \cite{Pasquale1999} and the very recent paper \cite{Schlichtkrull-Trapa-Vogan2018}. 
\vskip 8pt
\noindent Second, the degenerate principal series representations of the complex symplectic Lie group $\spg(n,\C)$ are "small" in the sense of the Kirillov-Gelfand dimension(see e.g. \cite{Vogan2017}) amongst infinite-dimensional unitary representations of $\spg(n,\C)$. According to the guiding  principle "small representation of a group $=$ large symmetries in a representation space" suggested by T. Kobayashi in \cite{Kobayashi2013}, explicit models of such representations are a natural source of information on special functions that arise in this framework as specific vectors. This philosophy has been applied to the analysis of minimal representations of $\org(p,q)$ (see \cite{Kobayashi-Mano2011} and \cite{Kobayashi-Orsted2003}) and small principal series representations of the real symplectic group $\spg(n,\R)$ (see \cite{Kobayashi-Orsted-Pevzner2011}).
\vskip 8pt
\noindent Fix an integer $n \geq 2$, set $m=n-1$ and write $G=\spg(n,\C)$ and $K=\spg(n)$. The representations we work with are defined by parabolic induction with respect to a maximal parabolic subgroup $Q$ of $G$ whose Langlands decomposition is $Q=MAN$(see Section \ref{subsection-non_comp._pict._and_heisenberg_group} for explicit description) with:
$$M \simeq \ung(1) \times \spg(m,\C), \; A \simeq \R^{\times}_{+} \; {\rm and} \; N \simeq {\rm H}_{\C}^{2m+1}.$$
Here, $\mathrm{H}_{\C}^{2m+1}$ refers to the $(2m+1)$-dimensional complex Heisenberg group. We denote the $\ung(1)$ component of an element $m$ of $M$ by $e^{i \theta(m)}$ and the positive real scalar that corresponds to an element $a$ of $A$ by $\alpha(a)$. Now consider for $(\lambda,\delta) \in \R \times \Z$ the character $\chi_{i\lambda,\delta}$ defined on $Q$ by:
$$\chi_{i\lambda,\delta}(man)= \big( e^{i \theta(m)}\big)^{\delta}  \,  \big( \alpha(a)\big)^{i\lambda}.$$
The corresponding induced representation $\pild={\rm Ind}_Q^G \, \chi_{i\lambda,\delta}$ can be realised on the completion $V_{i\lambda,\delta}$ of the complex vector space
{\small $$V_{i\lambda,\delta}^0 = \left \lbrace f \in C^0 \big( \C^{2n} \setminus \{0\} \big) \; \Big/ \; \forall c \in \C \setminus \{0\}: \; f(c \, \cdot)=\left( \frac{c}{\vert c \vert} \right)^{-\delta} \vert c \vert^{-i\lambda-2n} f(\cdot) \right \rbrace$$}with respect to the $L^2$-norm on $S^{4n-1}$.
\vskip 8pt
\noindent Representations $\pild$ form a degenerate principal series of $G$. It is proved in \cite{Gross1971} that $\pild$ is irreducible if $(\lambda,\delta) \neq (0,0)$ and in \cite{Clare2012} that $\pi_{0,0}$ decomposes into a sum of two irreducible subrepresentations. The isotypic decomposition of $\pild \big|_K$ is multiplicity free and given by (see \cite{Clare2012} and \cite{Howe-Tan1993}):
\begin{equation} \label{equation-general_isot._dec.}
\pild \big|_K \simeq \sideset{}{^\oplus}\sum_{ \substack{l-l' \geq |\delta| \\
l-l' \equiv \delta [2]}} \,  \pi^{l,l'}.
\end{equation}
In this sum, $l$ and $l'$ are integers such that $l \geq l' \geq 0$ and $\pi^{l,l'}$ is the irreducible representation of $K$ whose highest weight is $(l,l',0,\cdots,0)$; we denote the irreducible invariant subspace of $V_{i\lambda,\delta}$ (with respect to the left action of $K$) that corresponds to $\pi^{l,l'}$ by $V^{l,l'}$, calling it a \textit{component} of $V_{i\lambda,\delta}$.
\vskip 8pt
\noindent Our aim is to describe specific elements of components $V^{l,l'}$ in terms of special functions. This will require suitable changes of the carrying space $V_{i\lambda,\delta}$ (together with the action of $G$) so as to have a clearer view of $\pild$, depending on the kind of special functions we have in mind; each point of view is called a \textit{picture} of $\pild$ (above definition is the \textit{induced picture}).
\vskip 8pt
\noindent This paper is organised as follows:
\begin{itemize}
\item Section 2: in the compact picture, we study the $K$-type structure, by which we mean the detailed description of the components of $V_{i\lambda,\delta}$ and connections that exist between them (Propositions  \ref{proposition-left_action_decomposition}, \ref{proposition-hwv_for_right_action} and \ref{proposition-isotypic_decomposition_for_sp(1)}).
\item Section 3: we consider the case $l=l'$ and use, again in the compact picture, invariance properties with respect to $\spg(1)$ and $1 \times \spg(n-1)$ to exhibit in components $V^{l,l}$ elements which can be seen as solutions of hypergeometric differential equations (Theorem \ref{theorem-hypergeometric_equation}).
\item Section 4: we define the non-standard picture of $\pild$ (which was introduced in \cite{Kobayashi-Orsted-Pevzner2011} and followed in  \cite{Clare2012}) by applying a certain partial Fourier transform $\mathcal{F}$ to the non-compact picture. For a wide class of components, namely components $V^{l,0}$ (said otherwise, those components such that $l'=0$), this leads to elements that can be expressed in terms of modified Bessel functions (Theorem \ref{theorem-final_formula}).
\end{itemize}
Let us put together our most important results:
\begin{nnthm}[Main results] \label{thm-quat._sph._harm._thm} \

\noindent Let $n \in \N$ be such that $n \geq 2$ and set $m=n-1$.\\
\noindent Consider the group $G=\spg(n,\C)$, its maximal compact subgroup $\spg(n)$ and the parabolic subgroup $Q=MAN \simeq \ung(1) \times \spg(m,\C) \times \R^{\times}_{+} \times {\rm H}_{\C}^{2m+1}$.\\
\noindent Consider a pair $(\lambda,\delta) \in \R \times \Z$, together with the character $\chi_{i\lambda,\delta}$ defined by $\chi_{i\lambda,\delta}(man)= \big( e^{i \theta(m)}\big)^{\delta}  \,  \big( \alpha(a)\big)^{i\lambda}$
and the degenerate principal series representations $\pi_{i\lambda,\delta} = {\rm Ind}_Q^G \, \chi_{i\lambda,\delta}$ of $G$. 
\begin{enumerate}
\item For $l \in \N$, consider the $K$-type $\pi^{l,l}$. The corresponding subspace $V^{l,l}$ of $V_{i\lambda,\delta}$ contains an element which can be seen as a function of a single variable $\tau \in [0,1]$ and the restriction $\varphi$ of this function to $]0,1[$ satisfies the following hypergeometric equation:
$$\quad \; \tau(1-\tau) \varphi''(\tau) \, + \, 2(1-n\tau) \, \varphi'(\tau) \, + \, l(l+2n-1) \, \varphi(\tau) \, = \, 0.$$
\item For $(l,\alpha,\beta) \in \N^3$ such that $l=\alpha+\beta$ and $\delta=\beta-\alpha$, consider the $K$-type $\pi^{l,0}$ of $\pild$. Then, in the non-standard picture, highest weight vectors of $\pi^{l,0}$ are proportionnal to a function
$$\psi \, : \, \C \times \C^m \times \C^m \, \longrightarrow  \C$$
whose expression for $s \neq 0$ and $v \neq 0$ is
$$\quad \psi(s,u,v) \, = \, R(s,u,v) \, K_{\frac{i\lambda+\delta}{2}} \left( \pi \sqrt{ 1 + \Vert u \Vert^2 } \sqrt{ \vert s \vert^2+4\Vert v \Vert^2 } \right)$$
where we set
$$\; R(s,u,v) \, = \, \frac{(-i \overline{s})^{\alpha} \, \pi^{i\lambda+\beta+n}}{2^{\frac{i\lambda+l}{2}+1} \, \Gamma \left( \frac{i\lambda+l}{2}+n \right)} \, \left( \frac{ \sqrt{ |s|^2+4\Vert v \Vert^2}}{\pi \sqrt{1+\Vert u \Vert^2}} \right)^{\frac{i\lambda+\delta}{2}}$$
and where $K_{\frac{i\lambda+\delta}{2}}$ denotes a modified Bessel function of the third kind (see appendix for definition).
\end{enumerate}
\end{nnthm}
\noindent Before we enter the details, let us motivate the use of partial Fourier transforms. For one thing, they have proved useful in the study of Knapp-Stein operators (see \cite{Clare2012}, \cite{Kobayashi-Orsted-Pevzner2011}, \cite{Pevzner-Unterberger2007} and \cite{Unterberger2003}). For another, partial Fourier  transforms with respect to appropriate Lagrangian subspaces modify the nature of constraints imposed on specific vectors in representation spaces and have been used to find explicit formulas for $K$-finite vectors in \cite{Kobayashi-Mano2011}, \cite{Kobayashi-Orsted2003},  and \cite{Kobayashi-Orsted-Pevzner2011}. These works  establish formulas that involve Bessel functions. This has lead us to apply similar techniques to degenerate principal series of $\spg(n,\C)$.
\section{$K$-type structure} \label{section-k_type_structure}
\noindent We investiqate the $K$-type structure of the degenerate principal series $\pild={\rm Ind}_Q^G \, \chi_{i\lambda,\delta}$ of $\spg(n,\C)$. 
\subsection{Compact picture and general facts} \label{subsection-compact_picture_and_general_facts} 
\
\vskip 8pt
\noindent Let us fix $(\lambda,\delta) \in \R \times \Z$. In the compact picture (see \cite{Knapp1986}, Chapter VII), the carrying space of $\pild$ is a subspace of $L^2(G)$. The natural action of $G$ on $\C^{2n}$ enables one to identify it with the following Hilbert space:
$$L^2_{\delta}(S^{4n-1})=\left \lbrace f \in L^2(S^{4n-1}) \; \big/ \; \forall \theta \in \R: \; f(e^{i\theta} \, \cdot)=e^{-i\delta \theta}  f(\cdot) \right \rbrace.$$
\vskip 8pt
\noindent The compact picture is the ideal setting to look for irreducible invariant subspaces with respect to the left action of $K$, because one can benefit from decompositions given by standard harmonic analysis (see e.g. Chapter 9 of \cite{Faraut2008}, Chapters IV and V of \cite{Knapp2002} and Chapters 9 and 11 of \cite{Vilenkin1993-2}).
\vskip 8pt
\noindent Denote by $\harm^k$ the complex vector space of polynomial functions $f$ on $\R^{4n}$ that are  harmonic and homogeneous of degree $k$. Let us write $\sph^k=\harm^k \big|_{S^{4n-1}}$ (elements of $\sph^k$ are called spherical harmonics). Then:
\begin{itemize}
\item each $\sph^k$ is invariant under the left action of $\sog(4n)$ (this is also true for each $\harm^k$);
\item the representations defined by the left action of $\sog(4n)$ on the various spaces $\sph^k$ are irreducible and pairwise inequivalent;
\item $\displaystyle L^2\left( S^{4n-1} \right) = \widehat{\bigoplus_{k \in \N}} \; \sph^k$.
\end{itemize}
Consider the identification $(x,y) \in \R^{2n}\times \R^{2n} \longleftrightarrow z=x+iy \in \C^{2n}$. Then:
\begin{itemize}
\item functions of the variables $x$ and $y$ (in particular elements of the spaces $\harm^k$) can be regarded as functions of the variables $z$ and $\bar{z}$;  
\item matrices $A+iB$ of $\glg(2n,\C)$ (resp. $\sug(2n)$), where $A$ and $B$ denote $2n \times 2n$ real matrices, can be regarded as matrices $\left( \begin{array}{cc}
       A & -B \\
       B & A \\
       \end{array}
\right)$ of $\glg(4n,\R)$ (resp. $\sog(4n)$);
\item accordingly, the left action $L$ of $\sug(2n)$ on functions of $z$ and $\bar{z}$ is defined by $L(u)f(z,\bar{z}) = f(u^{-1}z, \, \overline{u^{-1} z}) = f(u^{-1}z,\, {^t u}\bar{z})$.
\end{itemize}
\vskip 8pt
\noindent The Laplace operator can be written $\displaystyle \Delta =4\sum_{i=1}^{2n} \frac{\partial^2}{\partial z_i \partial \bar{z_i}}$. For $\alpha,\beta \in \N$, consider the space $\harm^{\alpha,\beta}$
of polynomial functions $f$ of the variables $z$ and $\bar{z}$ such that $f$ is homogeneous of degree $\alpha$ in $z$ and degree $\beta$ in $\bar{z}$ and such that $\Delta f=0$. Let us write $\sph^{\alpha,\beta}=\harm^{\alpha,\beta} \big|_{S^{4n-1}}$. Then:
\begin{itemize}
\item each $\sph^{\alpha,\beta}$ is invariant under the left action of $\sug(2n)$ (this is also true for each $\harm^{\alpha,\beta}$);
\item the representations defined by the left action of $\sug(2n)$ on the various spaces $\sph^{\alpha,\beta}$ are irreducible and pairwise inequivalent;
\item $\displaystyle \sph^k = \bigoplus_{
\substack{
(\alpha,\beta) \in \N^2 \\
\alpha+\beta=k}} \sph^{\alpha,\beta}$.
\end{itemize}
This leads to the Hilbert sum $\displaystyle L^2_{\delta}(S^{4n-1})=\widehat{\bigoplus_{
\substack{(\alpha,\beta) \in \N^2 \\
\delta=\beta-\alpha}}}  \sph^{\alpha,\beta}$. We see that, in order to describe the isotypic decomposition of $\pild \big|_K$, we need to understand how each $\sph^{\alpha,\beta}$ breaks into irreducible invariant subspaces under the left action of $K$.
\vskip 8pt
\noindent From now on, we consider $3$-tuples $(k,\alpha,\beta) \in \N^3$ such that $\alpha+\beta=k$ and denote by $L$ the left action of $K$, be it on $\harm^k$, $\sph^k$, $\harm^{\alpha,\beta}$ or $\sph^{\alpha,\beta}$.
\subsection{Left action of $\spg(n)$} \label{subsection-left_action}
\
\vskip 8pt
\noindent Recall that the Lie algebra $\spa(n,\C)$ of $G$ is
$$\mathfrak{g}=\left\lbrace X=\left(
\begin{array}{cc}
A & C \\
B & -\trans{A} \\
\end{array}
\right) \in \m(2n,\C)\ / \ B\ {\rm and}\ C\ {\rm are\  symmetric} \right\rbrace$$ 
and that the Lie algebra $\spa(n)$ of $K$ is
$$\mathfrak{k}=\left\lbrace X=\left(
\begin{array}{cc}
A & -\overline{B} \\
B & \overline{A} \\
\end{array}
\right) \in \m(2n,\C)\ / \ A\ {\rm is\ skew \; and}\ B\ {\rm is\  symmetric} \right\rbrace.$$ 
\noindent Consider the complexification $\mathfrak{g}=\mathfrak{k}\oplus i\mathfrak{k}$ of $\mathfrak{k}$. Let $\mathfrak{h}$ be the usual Cartan subalgebra of $\mathfrak{g}$ consisting of diagonal elements of $\mathfrak{g}$. If $r$ belongs to $\{1,...,n\}$, denote by $L_r$ the linear form that assigns to an element of $\mathfrak{h}$ its $r^{\rm th}$ diagonal term.
\vskip 8pt
\noindent The set $\Delta$ of roots of $\mathfrak{g}$ (with respect to $\mathfrak{h}$) consists of the following linear forms:
\begin{itemize}
\item $L_r-L_s$ and $-L_r+L_s$ ($1 \leq r <s \leq n$);
\item $L_r+L_s$ and $-L_r-L_s$ ($1 \leq r <s \leq n$);
\item $2L_r$ and $-2L_r$ ($1 \leq r \leq n$).
\end{itemize}
\noindent The corresponding root spaces are respectively generated by:
\begin{itemize}
\item $U^+_{r,s}=E_{r,s} - E_{n+s,n+r}$ and $U^-_{r,s}=E_{s,r} - E_{n+r,n+s}$;
\item $V^+_{r,s}=E_{r,n+s} + E_{s,n+r}$ and $V^-_{r,s}=E_{n+s,r} + E_{n+r,s}$;
\item $D_r^+=E_{r,n+r}$ and $D_r^-=E_{n+r,r}$.
\end{itemize}
\noindent Here, $E_{r,s}$ denotes the elementary $2n \times 2n$ matrix whose coefficients are all $0$ except for the one that sits in line $r$ and column $s$ and which is taken to be $1$. The set of positive roots is
$$\Delta^+ \, = \, \lbrace L_r-L_s,L_r+L_s \rbrace_{1 \leq r <s \leq n} \cup \lbrace 2L_r\rbrace_{1 \leq r \leq n}$$
and the corresponding root spaces are those generated by the various matrices $U^+_{r,s}$, $V^+_{r,s}$ and $D_r^+$.
\begin{nnrem} Restriction to the unit sphere $S^{4n-1}$ intertwines the left action of $\sog(4n)$ on $\harm^k$ (resp. of $\sug(2n)$ on $\harm^{\alpha,\beta}$) and the left action of $\sog(4n)$ on $\sph^k$ (resp. of $\sug(2n)$ on $\sph^{\alpha,\beta}$). It is therefore equivalent for us to study the left action $L$ of $K$ on polynomials of $\harm^k$ (resp. $\harm^{\alpha,\beta}$) or on functions of $\sph^k$ (resp. $\sph^{\alpha,\beta}$); for convenience, we will work with polynomials.
\end{nnrem}
\noindent The infinitesimal action of $K$ is defined for  $P \in \harm^{\alpha,\beta}$ and $X \in \mathfrak{k}$ by:
$$dL(X)(P)(z,\bar{z})=\frac{d}{dt}\Big|_{t=0} \big( P\left( \exp(-tX) z, \trans{ \big( \exp(tX) \big) } \bar{z} \right) \big).$$
This formula can be written:
\begin{equation} \label{equation-infinitesimal_action}
dL(X)(P)(z,\bar{z})=
\left( \begin{array}{cc}
       \frac{\partial P}{\partial z}(z,\bar{z}) & \frac{\partial P}{\partial \bar{z}}(z,\bar{z})
       \end{array}
\right)
\left( \begin{array}{c}
       -Xz \\
       {^t X} \bar{z}
       \end{array}
\right).
\end{equation}
Extend the infinitesimal action $dL$ to $\mathfrak{g}$: this complexification is also defined by Formula (\ref{equation-infinitesimal_action}), but this time with $X \in \mathfrak{g}$.
\vskip 8pt
\noindent To make polynomial calculations clearer, from now on we change our system of notation for complex variables. Throughout the rest of this paper, we denote by $(z_1,...,z_n,w_1,...,w_n)$ the coordinates of $\C^{2n}$, or simply by $(z,w)$, in which case it is understood that $z=(z_1,...,z_n)$ and $w=(w_1,...,w_n)$. 
\begin{prop} \label{proposition-left_action_decomposition}
Let $(k,\alpha,\beta) \in \N^3$ be such that $k=\alpha+\beta$ and consider the left action of $K$ on $\harm^{\alpha,\beta}$. Given any integer $\gamma \in [0,{\rm min}(\alpha,\beta)]$, the polynomial 
$$P_{\gamma}^{\alpha,\beta}(z,w,\bar{z},\bar{w})=w_1^{\alpha - \gamma}\bar{z_1}^{\beta - \gamma}(w_2\bar{z_1}-w_1\bar{z_2})^{\gamma}$$
is a highest weight vector whose highest weight is the linear form
$$\sigma_{\gamma}^k = (k-\gamma) L_1+\gamma L_2.$$
Refering to the basis $\lbrace L_1,\cdots,L_n \rbrace$ of $\mathfrak{h}^\ast$, we will write:
$$\sigma_{\gamma}^k=(k-\gamma,\gamma,0,\cdots,0).$$
\noindent We shall denote by $V_{\gamma}^{\alpha,\beta}$ the corresponding irreducible invariant subspace of $\harm^{\alpha,\beta}$. One then has:
$$\displaystyle  \harm^{\alpha,\beta}=\bigoplus_{\gamma =0}^{{\rm min}(\alpha,\beta)}\,  V_{\gamma}^{\alpha,\beta}.$$
\end{prop}
\begin{flushleft}
\textbf{Proof: }
\end{flushleft}
Straightforward calculations, using (\ref{equation-infinitesimal_action}), show that $P_{\gamma}^{\alpha,\beta}$ is indeed a highest weight vector whose highest weight is $\sigma_{\gamma}^k$. Applying Weyl's dimension formula to the highest weight $\sigma_{\gamma}^k$, one can check the standard formula for the dimension $d_{\gamma}^k$ of $V_{\gamma}^{\alpha,\beta}$:
$$d_{\gamma}^k = \frac{(k-\gamma+2n-2)! \,(\gamma+2n-3)!\, (k-2\gamma+1) \,(k+2n-1)}{(k-\gamma+1)! \, \gamma! \, (2n-1)! \,(2n-3)!}.$$
One can then show by induction that the sum of these dimensions matches the dimension of $\harm^{\alpha,\beta}$, which finishes the proof (due to the Peter-Weyl theorem).
\begin{flushright}
\textbf{End of proof.}
\end{flushright}
\begin{nnrem} we have re-established the isotypic decomposition (\ref{equation-general_isot._dec.}), by identifying the highest weights $(l,l',0,\cdots,0)$ of (\ref{equation-general_isot._dec.}) with the highest weights $(k-\gamma,\gamma,0,\cdots,0)$ of Proposition \ref{proposition-left_action_decomposition}.
\end{nnrem}
\subsection{Right action of $\spg(1)$} \label{subsection-right_action}
\
\vskip 8pt
\noindent An interesting relation exists between the  various $K$-types, one that corresponds to a scalar type multiplication by elements of $\spg(1)$ which is inspired by the following decompositions:
\begin{align*}
L^2(S^{4n-1}) & =L^2_{\rm even}(S^{4n-1}) \oplus L^2_{\rm odd}(S^{4n-1}) \\
L^2(S^{4n-1}) & = \widehat{\bigoplus_{\delta \in \Z}} \; \; L^2_{\delta}(S^{4n-1}).
\end{align*}
Above summands are invariant under scalar multiplication of coordinates by unit numbers of respectively $\R$ and $\C$; the parameters even/odd and $\delta$ correspond to the characters of respectively $\org(1)$ and $\ung(1)$. So two actions seem to rule the decomposition of $L^2(S^{4n-1})$: left action of matrices and scalar multiplication by unit numbers. This diagram captures the situation:
$$\begin{array}{ccccc}
\org(4n) & \curvearrowright & L^2(S^{4n-1}) &  \curvearrowleft & \org(1) \\
\bigcup & & & &  \bigcap \\
\ung(2n) & \curvearrowright & L^2(S^{4n-1}) &  \curvearrowleft & \ung(1).
\end{array}$$
\noindent One would like to add the line
$$\begin{array}{ccccc}
\spg(n) & \curvearrowright & L^2(S^{4n-1}) &  \curvearrowleft & \spg(1).
\end{array}$$
\noindent This suggests an interaction between the left action of $\spg(n)$ and a scalar-type action of the group $\spg(1)$.
\vskip 8pt
\noindent This action of $\spg(1)$ is best understood in the quaternionic setting. In this work, quaternions $q$ are written $q=a+jb$, with $(a,b) \in \C^2$ and $j^2=-1$. The modulus of a quaternion $q$ is $|q|=\sqrt{|a|^2+|b|^2}$ and we say that $q$ is a unit quaternion if $|q|=1$. The field $\quat$ of quaternions is not abelian, because of the multiplication rule $jz=\overline{z}j$ ($z \in \C)$.
\vskip 8pt
\noindent The vector space $\quat^n$ can be seen as a left or right vector space; we choose to regard it as a right one. We find this choice convenient because applying matrices on the left of  column vectors commutes with scalar multiplication of coordinates on the right. Vectors $h \in \quat^n$ will be written $h=z+jw$, with $(z,w) \in \C^n \times \C^n$; the norm of such a vector $h$ is then $\Vert h \Vert=\sqrt{\Vert z \Vert^2 + \Vert w \Vert^2}$.
\vskip 8pt
\noindent Objects in the quaternionic setting can be read in the complex setting:
\begin{itemize}
\item functions of the variable $h \in \quat^n$ can be regarded as functions of the variables $z$ and $w$;  
\item matrices $A+jB$ of $\glg(n,\quat)$, where $A$ and $B$ denote $n \times n$ complex matrices, can be regarded as matrices $\left( \begin{array}{cc}
       A & -\overline{B} \\
       B & \overline{A} \\
       \end{array}
\right)$ of $\glg(2n,\C)$. In particular:
\begin{itemize}
\item the group of linear isometries of $\quat^n$ identifies with $\spg(n)$;
\item the group of unit quaternions $q=a+jb$ identifies with the group of $2 \times 2$ complex matrices $\left( \begin{array}{cc}
       a & -\overline{b} \\
       b & \overline{a} \\
       \end{array}
\right)$ such that $|a|^2+|b|^2=1$, that is, the group $\spg(1)$.
\end{itemize}
\end{itemize}
\vskip 8pt
\noindent The rules of quaternionic multiplication imply, given $h=z+jw \in \quat^n$ and $q=a+jb \in \quat$:
$$hq=(az-b\overline{w}) + j (aw+b\overline{z}).$$
One can transfer this multiplication to the complex setting. This is what motivates the way we define the right action of $\spg(1)$ on functions:
\begin{nndef}[Right action of $\spg(1)$] Consider a subset $\mathcal{S}$ of $\C^{2n}$ which is stable under above multiplication. Write elements of $\mathcal{S}$ as pairs $(z,w) \in \C^n \times \C^n$. Consider the set $\mathcal{F} = \lbrace f : \mathcal{S} \longrightarrow \C \rbrace$. Then the \textit{right action} $R$ of $\spg(1)$ on $\mathcal{F}$ is defined by:
$$R(q)f(z,w)=f \left(az-b\overline{w},aw+b\overline{z} \right)$$
for all $\big( q,f,(z,w) \big) \in \spg(1) \times \mathcal{F} \times \mathcal{S}$.
\end{nndef}
\begin{nnrem} Action $R$ preserves the space of homogeneous polynomials of degree $k$; what follows shows that it preserves the subspace $\harm^k$.
\end{nnrem}
\noindent On the Lie algebra level, by definition:
\begin{equation} \label{equation-right_action-before_extension}
dR(X)P(z,w)=\frac{d}{dt}\Big|_{t=0} \big( R \left( \exp(tX)\right) P(z,w) \big)
\end{equation}
where we take $P \in \harm^k$, $X \in \spa(1)$ and $(z,w) \in \C^n \times \C^n$. One extends $dR$ to the complexification $\sla(2,\C)$ of $\spa(1)$ by setting for $Z \in \sla(2,\C)$:
\begin{equation} \label{equation-right_action-extension}dR(Z)=dR(X)\, + \, i \, dR(Y)
\end{equation}
with $X$ and $Y$ in $\spa(1)$ and $Z=X+iY$. Let us consider the complex basis of $\sla(2,\C)$ that consists of the following matrices:
$$H=\left( 
\begin{array}{cc}
1 & 0 \\
0 & -1 \\
\end{array}
\right)\ \ \ \ \ \ 
E=\left( 
\begin{array}{cc}
0 & 1 \\
0 & 0 \\
\end{array}
\right)\ \ \ \ \ \ 
F=\left( 
\begin{array}{cc}
0 & 0 \\
1 & 0 \\
\end{array}
\right).$$
\begin{prop} \label{proposition-hwv_for_right_action} For $(\alpha,\beta,\gamma) \in \N^3$ such that $\gamma \leq {\rm min}(\alpha,\beta)$:
\begin{itemize} 
\item $dR(H)P_{\gamma}^{\alpha,\beta}=(\alpha-\beta) \, P_{\gamma}^{\alpha,\beta}.$ 
\item $dR(E)P_{\gamma}^{\alpha,\beta} = 
\left\{ \begin{array}{ll}
(\beta - \gamma) \, P_{\gamma}^{\alpha+1,\beta-1} & {\rm if} \, \gamma < \beta, \\
0 & {\rm if} \, \gamma = \beta.
\end{array} \right.$
\item $dR(F)P_{\gamma}^{\alpha,\beta} = 
\left\{ \begin{array}{ll}
(\alpha - \gamma) \, P_{\gamma}^{\alpha-1,\beta+1} & {\rm if} \, \gamma < \alpha, \\
0 & {\rm if} \, \gamma = \alpha.
\end{array} \right.$
\end{itemize}
Consequently, given $(k,\gamma) \in \N^2$ such that $\gamma \leq \E \left( \frac{k}{2} \right)$ (denoting by $\E$ the function that assigns to a real number its integer part), $P_{\gamma}^{k-\gamma,\gamma}$ is a highest weight vector of $R$ whose highest weight is $k-2\gamma$.
\end{prop}
\begin{flushleft}
\textbf{Proof: }
\end{flushleft}
Define the following elements of $\spa(1)$:
$$X=\left( 
\begin{array}{cc}
i & 0 \\
0 & -i \\
\end{array}
\right)\ \ \ \ \ \ 
Y=\left( 
\begin{array}{cc}
0 & i \\
i & 0 \\
\end{array}
\right)
\ \ \ \ \ \ 
Z=\left( 
\begin{array}{cc}
0 & -1 \\
1 & 0 \\
\end{array}
\right).$$
\noindent 
Write:
$$H = 0+i(-X)\ \ \ \ \ \ 
E = \frac{-Z}{2} +  i\left( \frac{-Y}{2} \right)\ \ \ \ \ \ 
F = \frac{Z}{2} +  i\left( \frac{-Y}{2} \right).$$
One then applies (\ref{equation-right_action-before_extension}) and (\ref{equation-right_action-extension}) to obtain the desired formulas.
\begin{flushright}
\textbf{End of proof.}
\end{flushright}
\subsection{K-type diagram} \label{subsection-k_type_diagram}
\begin{figure}
\centering
\includegraphics[scale=0.6]{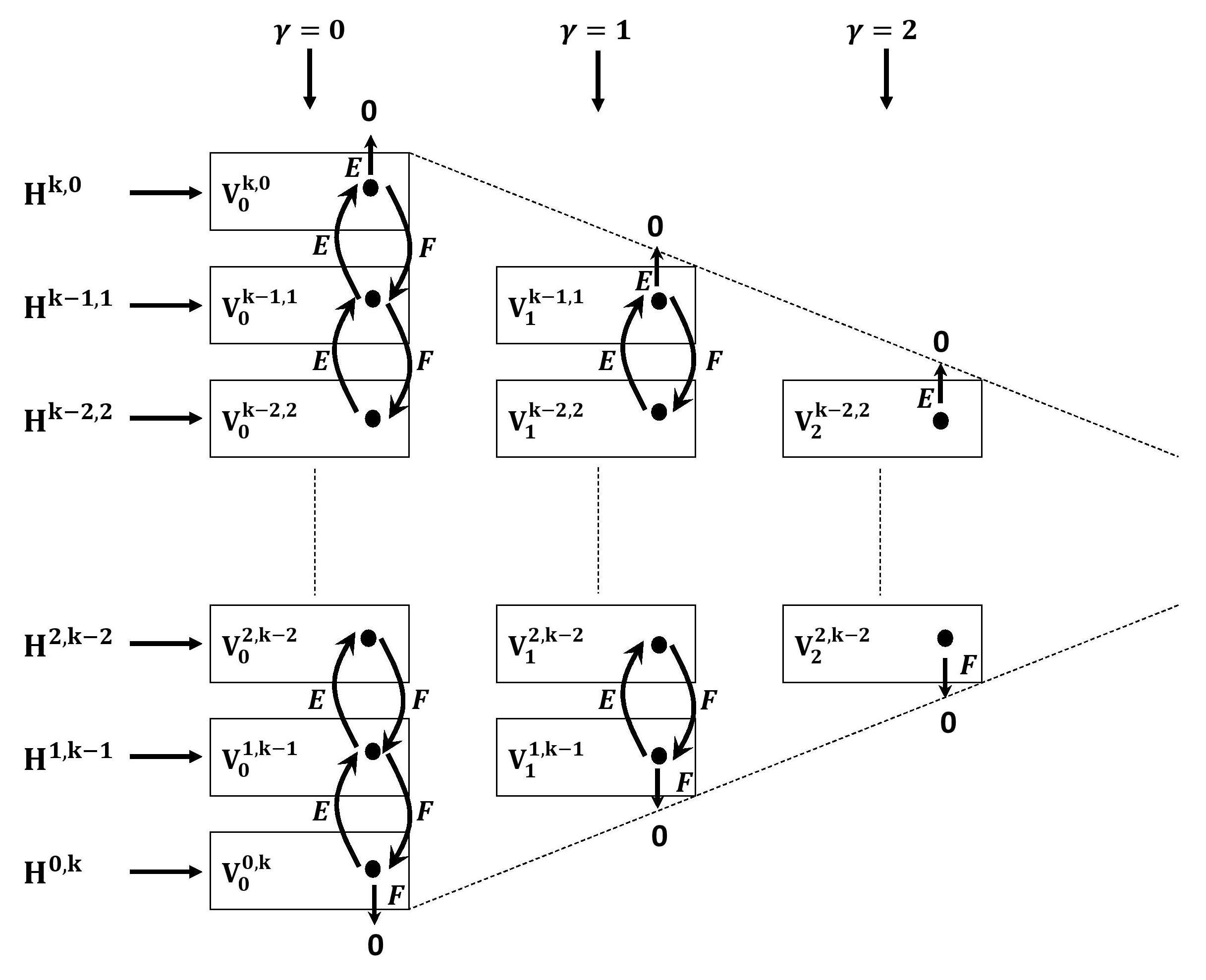}
\caption{$K$-type diagram}
\label{diagram-both_actions}
\end{figure}

\
\vskip 8pt
\noindent Figure \ref{diagram-both_actions} captures the contents of Propositions \ref{proposition-left_action_decomposition} and \ref{proposition-hwv_for_right_action}. In this diagram, the thick black dots represent highest weight vectors $P_{\gamma}^{\alpha,\beta}$ of $L$ (again, $\alpha+\beta = k$). Let us take a closer look at Figure \ref{diagram-both_actions}. For any fixed integer $\gamma$ such that $0 \leq \gamma \leq E \left( \frac{k}{2} \right)$:
\begin{itemize}
\item The arrow beneath a value of $\gamma$ points to a vertical set of components whose direct sum defines an invariant subspace $V_{\gamma}^k$ (obviously not irreducible) under the left action of $\spg(n)$:
$$V_{\gamma}^k = \bigoplus_{\alpha=\gamma}^{k-\gamma} V_{\gamma}^{\alpha,k-\alpha}.$$
\item A vertical set of thick black dots defines a basis of an irreducible invariant subspace $W_{\gamma}^k$ under the right action of $\spg(1)$:
$$W_{\gamma}^k = {\rm Vect}_{\C} \lbrace P_{\gamma}^{\alpha,k-\alpha} \rbrace_{\alpha = \gamma,\ldots,k-\gamma} \,\subset V_{\gamma}.$$
\item Because the left action of $\spg(n)$ and the right action of $\spg(1)$ commute, applying $L$ to $W_{\gamma}^k$ gives irreducible invariant subspaces of $R$ contained in $V_{\gamma}^k$; these subspaces define subrepresentations of $R$ which are equivalent to the restriction of $R$ to $W_{\gamma}^k$.
\end{itemize}
Using the fact that the left action $L$ of $\spg(n)$ is transitive in each component, we obtain:
\begin{prop} \label{proposition-isotypic_decomposition_for_sp(1)} The isotypic decomposition of the right action $R$ of $\spg(1)$ on the vector space $\harm^k$ is:
$$R \big|_{\harm^k} \simeq \sideset{}{^\oplus}\sum_{\quad 0 \leq \gamma \leq \E\left( \frac{k}{2} \right)} \, d_{\gamma}^k \, R \big|_{W_{\gamma}^k}.$$
We see that the multiplicity of $R \big|_{W_{\gamma}^k}$ is the dimension $d_{\gamma}^k$ of $V_{\gamma}^{k-\gamma,\gamma}$ (given in the proof of Proposition \ref{proposition-left_action_decomposition}). 
\end{prop}
\begin{nnrem}
The structure of $\harm^k$, which we have described with respect to the left action of $\spg(n)$ and the right action of $\spg(1)$, appears in \cite{Howe-Tan1993} (Proposition 5.1 and Lemma 6.1), where it is expressed in terms of tensor products:
$$\harm^k \big|_{\spg(n) \times \spg(1)}=\sum_{\gamma = 0}^{\E\left( \frac{k}{2} \right)} \, V_{\gamma}^{k-\gamma,\gamma} \otimes W_{\gamma}^k.$$
\end{nnrem}
\begin{nnrem}
Restriction to the unit sphere intertwines the right action of $\spg(1)$ on $\harm^k$ and the right action of $\spg(1)$ on $\sph^k$. Thus the structure of $\harm^k$ shown in Figure \ref{diagram-both_actions} is also the structure of $\sph^k$ under the right action of $\spg(1)$.
\end{nnrem}
\section{Compact picture and hypergeometric equations} \label{section-compact_picture_and_hyper._eq.}
\noindent The idea here is to consider eigenfunctions of the Casimir operator. We choose  eigenfunctions that have symmetry properties in order to reduce the number of variables.
\subsection{Bi-invariant functions} \label{subsection-bi_inv._functions}
\
\vskip 8pt
\noindent Proposition \ref{proposition-isotypic_decomposition_for_sp(1)} implies:
\begin{cor} \label{corollary-far_right_box_when_k_is_even}
Consider $k \in \N$. Existence in $\harm^k$ of a $1$-dimensional subspace that is stable under the right action of $\spg(1)$ implies that $k$ is even. Assuming now that $k$ is even and setting $\alpha = \frac{k}{2}$, $V_{\alpha}^{\alpha,\alpha}$ is the only component in $\harm^k$ to contain such a subspace; moreover, all elements of $V_{\alpha}^{\alpha,\alpha}$ are invariant under the right action of $\spg(1)$. 
\end{cor}
\noindent Let us now consider the subgroup $1 \times \spg(n-1)$ of $\spg(n)$ that consists of all matrices that can be written  
$\left( \begin{array} {cc}
1 & 0 \\
0 & A \\
\end{array} \right)$ in the quaternionic setting with $A \in \spg(n-1)$.
\begin{nndef}
Consider any $(\alpha,\beta) \in \N^2$. A polynomial of $\harm^{\alpha,\beta}$ and its corresponding spherical harmonic in $\sph^{\alpha,\beta}$ are said to be \textit{bi-invariant} if:
\begin{itemize}
\item they are invariant under the left action of $1 \times \spg(n-1)$;
\item they are also invariant under the right action of $\spg(1)$.
\end{itemize}
\end{nndef}
\begin{prop} \label{proposition-bi_invariant_polynomials} \
\begin{enumerate}
\item Consider any $(\alpha,\beta,\gamma) \in \N^3$ such that $\gamma \leq {\rm min}(\alpha,\beta)$.\\
\noindent Denote by ${\rm Inv}_{\gamma}^{\alpha,\beta}$ the complex vector space that consists of all polynomials of $V_{\gamma}^{\alpha,\beta}$ that are invariant under the left action of $1 \times {\rm Sp}(n-1)$. Then: 
$${\rm dim}_{\C}({\rm Inv}_{\gamma}^{\alpha,\beta})=\alpha+\beta-2\gamma+1.$$
We point out that this dimension is higher than or equal to $1$ and that it equals $1$ if and only if $\alpha=\beta=\gamma$.
\item Consider $k \in \N$. Existence in $\harm^k$ of a bi-invariant polynomial implies that $k$ is even. Assuming now that $k$ is even and setting $\alpha=\frac{k}{2}$, there exists in $\harm^k$ a unique (up to a constant) bi-invariant polynomial and this polynomial belongs to $V_{\alpha}^{\alpha,\alpha}$.
\end{enumerate}
\end{prop}
\begin{flushleft}
\textbf{Proof:}
\end{flushleft}
\noindent An invariant polynomial of $V_{\gamma}^{\alpha,\beta}$ defines a trivial $1$-dimensional representation of $1 \times {\rm Sp}(n-1)$, whose multiplicity in the left action of $1 \times {\rm Sp}(n-1)$ on $V_{\gamma}^{\alpha,\beta}$ is equal to ${\rm dim}_{\C}({\rm Inv}_{\gamma}^{\alpha,\beta})$. Using a classical branching theorem proved by Zhelobenko for $\spg(n)$ and $\spg(n-1)$ (Theorem 9.18 in \cite{Knapp2002}), one can compute this multiplicity and thereby prove Item (1). Item (2) just requires to combine Item (1) and Corollary \ref{corollary-far_right_box_when_k_is_even}.
\begin{flushright}
\textbf{End of proof.}
\end{flushright}
\vskip 8pt
\noindent Given a non-negative even integer $k$, in order to exhibit the unique (up to a constant) bi-invariant polynomial of $H^k$ given by Proposition \ref{proposition-bi_invariant_polynomials}, one can concentrate on polynomials of the following type (they are obviously bi-invariant and of total homogeneous degree $k$):
{\small $$P(z,w) = \sum_{A \in \mathcal{A}_k} \nu_A \, (z_1 \overline{z_1} + w_1 \overline{w_1})^a \, ( z_2 \overline{z_2} + \cdots + z_n \overline{z_n} + w_2 \overline{w_2} + \cdots w_n \overline{w_n})^b.
$$}where $\mathcal{A}_k= \lbrace (a,b) \in \N^2 \, / \, a+b=\frac{k}{2} \rbrace$ and where coefficients $\nu_A$ denote complex scalars (undetermined at first). We need  $P$ to be harmonic, which requires specific values for coefficients $\nu_A$; it is shown in \cite{Mendousse2017} (Section 3.2.3) how to compute them. 
\begin{nnexams} Let us write $U=z_1 \overline{z_1} + w_1 \overline{w_1}$ and $\displaystyle V=\sum_{r=2}^n z_r \overline{z_r} + w_r \overline{w_r}$. Then the unique bi-invariant polynomial (up to a contant) of $H^k$ is: 
\begin{itemize}
\item $(1-n) \, U \, + \, V$, when $k=2$;
\item $\frac{(2n-1)(n-1)}{3}  \, U^2 \, + \, (1-2n) \, UV \, + \, V^2$, when $k=4$.
\end{itemize}
\end{nnexams}
\subsection{Casimir operator of $\spg(n)$} \label{subsection-casimir_operator}
\
\vskip 8pt
\noindent We consider the inner-product defined on $\mathfrak{k}$ by:
$$(X,Y) \longmapsto -\frac{1}{2}\tra(XY).$$
where $\tra$ denotes the trace form. Taking $r$ and $s$ to denote integers that belong to $\{1,...,n\}$, define:
\begin{itemize}
\item $A_r = iE_{r,r} - iE_{n+r,n+r}$;
\item $B_{r,s} = E_{r,s} - E_{s,r} + E_{n+r,n+s} - E_{n+s,n+r}$ ($r \neq s$);
\item $C_{r,s} = iE_{r,s} + iE_{s,r} - iE_{n+r,n+s} - iE_{n+s,n+r}$ ($r \neq s$);
\item $D_r = E_{n+r,r} - E_{r,n+r}$;
\item $E_r = iE_{n+r,r} + iE_{r,n+r}$;
\item $F_{r,s}=E_{n+r,s} + E_{n+s,r} - E_{r,n+s} - E_{s,n+r}$ ($r \neq s$);
\item $G_{r,s}=iE_{n+r,s} + iE_{n+s,r} + iE_{r,n+s} + iE_{s,n+r}$ ($r \neq s$).
\end{itemize}
\noindent One can check that the set 
$$\mathcal{B}_{\mathfrak{k}} = 
\{ A_r,D_r,E_r \}_{r \in \{1,...,n\}} 
\cup 
\left\lbrace   \frac{B_{r,s}}{\sqrt{2}},\frac{C_{r,s}}{\sqrt{2}},\frac{F_{r,s}}{\sqrt{2}},\frac{G_{r,s}}{\sqrt{2}} \right\rbrace_{1 \leq r < s \leq n}$$
is an orthonormal basis of $\mathfrak{k}$. One can apply to elements of $\harm^k$ any sequence of operators $dL(X)$ (taking $X \in \mathfrak{k}$). In particular, let us consider the Casimir operator $\Omega_L$, which can be defined by:
\begin{equation} \label{equation-casimir_definition}
\Omega_L=\sum_{X \in \mathcal{B}_{\mathfrak{k}}} \big( dL(X) \big)^2.
\end{equation}
This operator commutes with $L$. So Schur's lemma implies that $\Omega_L$ is equal to a scalar multiple of the identity map when restricted to an irreducible invariant subspace. This scalar can be computed (see \cite{Mendousse2017}, Corolloray 3.3):
\begin{prop} \label{proposition-casimir_eigenvalue}
Consider the left action $L$ of $K$ on $V_{\gamma}^{\alpha,\beta}$. Denoting by $Id$ the identity map:
$$\Omega_L = - \big( (k - \gamma )^2 + 2nk + \gamma^2 - 2\gamma \big) \, Id.$$
\end{prop}
\subsection{From bi-invariance to hypergeometric equations} \label{subsection-from_bi_inv._to_hyper._eq.}
\
\vskip 8pt
\noindent In this section, we consider the quaternionic projective space $\qps$, whose elements we define with respect to right-multiplication: given $x \in \quat^n$, the \textit{quaternionic line} through $x$ is the subspace 
$$x\quat=\{ xh\ /\ h \in \quat \}.$$
\noindent Quaternionic matrices act naturally on quaternionic lines: an invertible matrix $M \in \glg(n,\quat)$ assigns to a line $x\quat$ the line $(Mx)\quat$.  We now study the orbits under the restriction of this action to $1 \times \spg(n-1)$. Given $x \quat \in {\rm P}^{n-1}(\quat)$, we denote by $\mathcal{O}(x \quat)$ the orbit of $x \quat$. We show in the next proposition that the orbits can be parametrised by a single real variable.
\begin{prop} \label{proposition-theta_variable}
Consider any $x = (x_1,...,x_n) \in \quat^n$. There exists a unique $\theta \in \left[ 0 , \frac{\pi}{2} \right]$ such that, writing $x(\theta)=(\cos\theta,0,...,0,\sin\theta)\in \mathbb{H}^n$, the quaternionic line $x(\theta) \quat$ belongs to $\mathcal{O}(x \quat)$. Moreover, $\theta$ is explicitly given by the following formulas ($\Vert \cdot \Vert$ denotes the norm of $\quat^{n-1}$):
\begin{itemize}
\item if $x_1 \neq 0$, then $\theta = \arctan \left( \Vert x' \Vert \right)$, taking $x'=\left( \frac{x_2}{x_1},...,\frac{x_n}{x_1} \right)$;
\item if $x_1=0$, then $\theta = \frac{\pi}{2}$.
\end{itemize}
\end{prop}
\begin{flushleft}
\textbf{Proof:}
\end{flushleft}
\noindent Suppose first that $x_1 \neq 0$. Then, writing $x' = \left( x_2 \frac{1 }{x_1},...,x_n \frac{1}{x_1} \right)$, we have $x \quat=(1,x')\mathbb{H}$. Because $\spg(n-1)$ acts transitively on spheres of $\quat^{n-1}$, the group $1 \times \spg(n-1)$ contains an element that takes $(1,x') \mathbb{H}$ onto $\left( 1,0,...,0,\Vert x' \Vert \right) \quat$. One then chooses $\theta \in \left[0 , \frac{\pi}{2} \right[$ such that:
$$\left( 1,0,...,0,\Vert x' \Vert \right) \quat = (\cos\theta,0,...,0,\sin\theta) \quat.$$
\noindent Suppose now that $x_1=0$. Then $\Vert x' \Vert = 1$ and one can choose an element of $1 \times \spg(n-1)$ that carries $x \quat$ onto $(0,...0,1)\quat$, leading to $\theta=\frac{\pi}{2}$. This establishes existence and formulas for the parameter $\theta$; proof of uniqueness is straightforward.
\begin{flushright}
\textbf{End of proof.}
\end{flushright}
\vskip 8pt
\noindent Let us assume, throughout this section that $k$ is even and set $\alpha = \frac{k}{2}$. Consider the space $\harm^k$ and its unique, up to a constant, bi-invariant polynomial $P$, which in fact belongs to $V_{\alpha}^{\alpha,\alpha} \subset \harm^{\alpha,\alpha}$ (see Proposition \ref{proposition-bi_invariant_polynomials}). Denote by $f_P$ the spherical harmonic that corresponds to $P$, that is:
$$f_P=P |_{S^{4n-1}} \, \in \, \sph^{\alpha,\alpha}.$$
\noindent Invariance under the right action of $\spg(1)$ enables $f_P$ to descend to a function $f$ on ${\rm P}^{n-1}(\quat)$:
$$\begin{array}{ccccl}
f & : & {\rm P}^{n-1}(\quat) & \longrightarrow & \C \\
  & & x \quat & \mapsto & f_P(x) {\rm ,} \, {\rm where}\, \, x \, \, {\rm is \, \, chosen \, \,  in} \, \, S^{4n-1}. \\
\end{array}
$$
\noindent Transferring, in the way one expects, the left action of $K$ on spherical harmonics to a left action, also denoted by $L$, of $K$ on complex-valued functions defined on ${\rm P}^{n-1}(\quat)$, we write: 
$$L(k)f (x \mathbb{H}) = f \big( (k^{-1} x) \mathbb{H} \big)$$
for all $(k,x) \in K \times S^{4n-1}$. Invariance of $f_P$ under the left action of $1 \times \spg(n-1)$ implies that $f$ descends to a function on the classifying set $\mathcal{O}$ of orbits, leading in turn to the following new function (via Proposition \ref{proposition-theta_variable}): 
$$\begin{array}{ccccl}
F & : & \left[ 0,\frac{\pi}{2} \right] & \longrightarrow & \C \\
 & & \theta & \mapsto & F(\theta)=f \big( x(\theta) \mathbb{H} \big) = f_P \big( x(\theta) \big).\\
\end{array}
$$
\vskip 8pt
\noindent The following theorem establishes the first part of our  "Main results" theorem (stated in the introduction):
\begin{thm} \label{theorem-hypergeometric_equation}
Fix any even non-negative integer $k$. Consider the unique (up to a constant) bi-invariant spherical harmonic $f_P$ of $\sph^k$ and the corresponding function $F$ (as introduced above).
\begin{itemize}
\item The restriction of $F$ to the open interval $\left] 0,\frac{\pi}{2} \right[$ satisfies the following hypergeometric differential equation:
$$\qquad \quad F''(\theta) \, + \, \left( \frac{6}{\tan(2 \theta)} + \frac{4n-8}{\tan \theta} \right) \, F'(\theta) \, + \, k(k+4n-2) \, F(\theta) \, = \, 0.$$
\item Consider the following smooth diffeomorphism (onto):
$$\begin{array}{ccccc}
\psi & : & \left] 0,\frac{\pi}{2} \right[ & \longrightarrow & ]0,1[ \\
 & & \theta & \longmapsto & \cos^2 \theta. \\
\end{array}$$
\noindent Then $\varphi=F \circ \psi^{-1}$ satisfies the following hypergeometric differential equation:
$$\qquad \quad \tau(1-\tau) \, \varphi''(\tau) \, + \, 2(1-n\tau) \, \varphi'(\tau) \, + \frac{k}{2}\left( \frac{k}{2}+2n-1 \right) \, \varphi(\tau) \, = \, 0.$$
\end{itemize}
\end{thm}
\begin{flushleft}
\textbf{Proof:}
\end{flushleft}
\noindent The Casimir operator $\Omega_L$ of $L$ is given by formula (\ref{equation-casimir_definition}). Given an element $x \in S^{4n-1}$, we have:
$$\Omega_L(f_P)(x) = \sum_{X \in \mathcal{B}_{\mathfrak{k}}} \frac{\partial^2}{\partial t^2}\Big|_{t=0} \Big( f_P \big( \exp(-tX) x \big) \Big).$$
\noindent Proposition \ref{proposition-casimir_eigenvalue} tells us that $V_{\alpha}^{\alpha,\alpha}$ is associated to the eigenvalue $\varpi=-\big( 2 \alpha^2+(4n-2) \alpha \big)$ of $\Omega_L$:
$$\Omega_L(f_P)(x) = \varpi \, f_P (x).$$
\noindent Given $\theta \in \left[ 0 , \frac{\pi}{2} \right]$, one can apply this equation to $x=x(\theta)$:
\begin{equation} \label{equation-casimir_equation-1}
\sum_{X \in \mathcal{B}_{\mathfrak{k}}} \frac{\partial^2}{\partial t^2}\Big|_{t=0} \Big( f_P \big( \exp(-tX) x(\theta) \big) \Big) = \varpi \, f_P \big( x(\theta) \big).
\end{equation}
\noindent One can convert this into a differential equation satisfied by the function $F$ of the real variable $\theta$. But this requires to know for each $X \in \mathcal{B}_{\mathfrak{k}}$ and each $t \in \R$, which parameter $\xi_{X,\theta}(t) \in \left[ 0, \frac{\pi}{2} \right]$ labels the orbit of $\exp(-tX) x(\theta) \mathbb{H}$ under the left action of $1 \times \spg(n-1)$. Proposition  \ref{proposition-theta_variable} gives us the value of $\xi_{X,\theta}(t)$. Then (\ref{equation-casimir_equation-1}) can be written
\begin{equation} \label{equation-casimir_equation-2}
\sum_{X \in \mathcal{B}_{\mathfrak{k}}} \frac{\partial^2}{\partial t^2} \Big|_{t=0} \Big( F \big( \xi_{X,\theta}(t) \big) \Big) = \varpi \, F(\theta).
\end{equation}
\noindent Because $f_P$ is invariant under the left action of $1 \times \spg(n-1)$, the only elements of $\mathcal{B}_{\mathfrak{k}}$ that actually contribute to  (\ref{equation-casimir_equation-2}) are (taking $2 \leq r \leq n$):
$$A_1 \ ,\ D_1 \ ,\ E_1 \ ,\ \frac{B_{1,r}}{\sqrt{2}} \ ,\ \frac{C_{1,r}}{\sqrt{2}} \ ,\ \frac{F_{1,r}}{\sqrt{2}} \ ,\ \frac{G_{1,r}}{\sqrt{2}}.$$
\noindent Denote  by $\mathcal{B}'_{\mathfrak{k}}$ the set of these matrices. Using the differentiation chain rule, (\ref{equation-casimir_equation-2}) becomes:
$$\sum_{X \in \mathcal{B}'_{\mathfrak{k}}} F'' \big( \xi_{X,\theta}(0) \big) \, \big( \xi_{X,\theta}'(0) \big)^2 + F' \big( \xi_{X,\theta}(0) \big) \, \xi_{X,\theta}^{''}(0) = \varpi \, F (\theta).$$
\noindent Finally, all we have to do now is determine the coordinates of the various $\exp(-tX) x(\theta)$ and deduce the terms $\xi_{X,\theta}(0)$, $\xi'_{X,\theta}(0)$ and $\xi''_{X,\theta}(0)$. Long-winded calculations (see \cite{Mendousse2017} for details) finally lead to the first hypergeometric differential equation; it is then straightforward to establish the second.
\begin{flushright}
\textbf{End of proof.}
\end{flushright}
\section{Non-standard picture and modified Bessel functions}\label{section-non_standard_picture_and_mod._bessel_functions}
\noindent As mentionned in the introduction, the non-standard picture relies on a certain partial Fourier transform, which we introduce below (following \cite{Kobayashi-Orsted-Pevzner2011} and \cite{Clare2012}). Such transforms reveal new aspects of functions, because they swap the roles of two fundamental operators: multiplication by coordinates and differentiation with respect to those coordinates. Also, integral formulas that result from applying partial Fourier transforms can lead to  classical integrals which may be computed in terms of special functions.
\subsection{Non-compact picture and complex Heisenberg group} \label{subsection-non_comp._pict._and_heisenberg_group}
\
\vskip 8pt
\noindent The general theory of induced representations, as presented in \cite{Knapp1986} (Chapter VII), includes, beside the induced and compact pictures of $\pild$, a so-called \textit{non-compact picture}. Its carrying space is $L^2(\overline{N})$, where, in our case, $\overline{N} = \{\trans{n} \, /\ n \in N \}$ ($\trans{n}$ denotes the transpose of $n$). We will realise this carrying space in the geometric setting of $\C^{2n}$. In order to do so, let us first describe our parabolic subgroup explicitely (we are aware that $m$ and $n$ refer to elements of subgroups as well as dimensions; context makes intended meanings clear):
\begin{itemize}
\item $M$ consists of matrices {\small $m=\left(
\begin{array}{cccc}
e^{i\theta(m)} & 0 & 0 & 0 \\
0 & A & 0 & C \\
0 & 0 & e^{-i\theta(m)} & 0 \\
0 & B & 0 & D \\
\end{array}
\right)$} such that $\theta (m) \in \R$ and $\left(
\begin{array}{cc}
A & C \\
B & D \\
\end{array}
\right) \in \spg(m,\C)$, where $A,B,C,D$ denote $m \times m$ complex matrices;
\item $A$ consists of matrices {\small $a=\left(
\begin{array}{cccc}
\alpha(a) & 0 & 0 & 0 \\
0 & I_m & 0 & 0 \\
0 & 0 & \big( \alpha(a) \big)^{-1} & 0 \\
0 & 0 & 0 & I_m \\
\end{array}
\right)$} such that $\alpha(a)>0$ ($I_m$ denotes the $m$-dimensional identity matrix);
\item $N$ consists of matrices {\small $n=\left(
\begin{array}{cccc}
1 & \trans{u} & 2s & \trans{v} \\
0 & I_m & v & 0 \\
0 & 0 & 1 & 0 \\
0 & 0 & -u & I_m \\
\end{array}
\right)$} such that $s\in \C$ and $u$ and $v$ both belong to $\C^m$.
\end{itemize}
Then $Q=MAN$ is indeed a parabolic subgroup of $G$.
\vskip 8pt
\noindent We regard the $(2m+1)$-dimensional complex Heisenberg group $\mathrm{H}_{\C}^{2m+1}$ as $\C \times \C^m \times \C^m$ equipped with the following product: 
$$(s,u,v)(s',u',v')=\left( s+s'+\frac{\langle v,u' \rangle - \langle u,v' \rangle}{2},u+u',v+v' \right).$$
Here, $\langle \cdot,\cdot \rangle$ is defined for elements $(x,y) \in \C^m \times \C^m$ by: 
$$\langle x,y \rangle = \sum_{i=1}^m x_i y_i.$$
We point out that, later on, this sum will be denoted by $x \cdot y$ instead of $\langle x,y \rangle$ when we consider $(x,y) \in \R^m \times \R^m$.
\vskip 8pt
\noindent We use the following sequence of identifications (with $n$ as above):
$$\begin{array}{ccccc}
\mathrm{H}_{\C}^{2m+1} & \simeq & \overline{N} & \simeq & {\rm H}\\
(s,u,v) & \longmapsto & \trans{n} &  \longmapsto & (1,u,2s,v). \\
\end{array}$$
where ${\rm H}=\{1\} \times \C^m \times \C \times \C^m$ is regarded as a complex hyperplane of $\C^{2n}$. The first identification map is a group isomorphism (one can also prove that $\mathrm{H}_{\C}^{2m+1} \simeq N$) and  the second identication results from the natural action of $\overline{N}$ on $\C^{2n}$ (applied to the element $(1,0,\cdots,0)$ of $\C^{2n}$).
\vskip 8pt
\noindent One can now identify $L^2(\overline{N})$ with  $L^2(\mathrm{H}_{\C}^{2m+1})$ and, given $f \in V_{i\lambda,\delta}^0$ (in the induced picture), regard the restriction $f|_{\rm H}$ as an element of $L^2(\mathrm{H}_{\C}^{2m+1})$.
\vskip 8pt
\noindent From now on, we write $\C^{2m+1}$ instead of $\mathrm{H}_{\C}^{2m+1}$.

\subsection{Definition of the non-standard picture} \label{subsection-definition_of_the_non_stand._pict.}
\
\vskip 8pt
\noindent We first define two partial Fourier transforms
\begin{align*}
\mathcal{F}_{\tau} & : L^2(\C^{2m+1}) \longrightarrow L^2(\C^{2m+1})\\
\mathcal{F}_{\xi} &: L^2(\C^{2m+1}) \longrightarrow L^2(\C^{2m+1})
\end{align*}
by setting for functions $g \in L^2(\C^{2m+1})$ that fulfill integrability conditions: 
$$\mathcal{F}_{\tau}(g)(s,u,v)=\int_{\C} g(\tau,u,v) \, e^{-2i\pi {\rm Re}(s\tau) } \, d\tau$$
$$\mathcal{F}_{\xi}(g)(s,u,v)=\int_{\C^m} g(s,u,\xi)\,e^{-2i\pi {\rm Re} \langle v,\xi \rangle}\,d\xi.$$
\noindent We then define the partial Fourier transform on which is based the non-standard picture:
$$\mathcal{F}= \mathcal{F}_{\tau} \circ \mathcal{F}_{\xi}.$$
\begin{nndef}
\noindent The \textit{non-standard picture} of $\pi_{i\lambda,\delta}$ has $L^2(\C^{2m+1})$ as carrying space. The action of $G$ is then the conjugate under $\mathcal{F}$ of the action of $G$ in the non-compact picture; in other words, $\mathcal{F}$ intertwines the action of $G$ in the non-compact picture and the action of $G$ in the non-standard picture.
\end{nndef}
\subsection{Connection with modified Bessel functions} \label{subsection-connection_with_modified_bessel_functions}
\subsubsection{Selected components} \label{subsubsection-selected_components}
\
\vskip 8pt
\noindent 
Because we intend to use the right action of $\spg(1)$, we consider an entire space $\harm^k$ of harmonic polynomials. We restrict our attention to those  components of Proposition \ref{proposition-left_action_decomposition} that are labelled by $\gamma=0$ and are subspaces of $\harm^k$, namely the components $V_0^{\alpha,\beta}$ with $\alpha + \beta=k$; the corresponding highest weight vectors are the polynomials 
$$P_0^{\alpha,\beta}(z,w,\bar{z},\bar{w})=w_1^{\alpha}\bar{z_1}^{\beta}.$$
\noindent We denote by $g_{\alpha,\beta}$ the restriction of $P_0^{\alpha,\beta}$ to the unit sphere $S^{4n-1}$. Let us call $g$ the function in the induced picture that corresponds to $g_{\alpha,\beta}$. It extends $g_{\alpha,\beta}$, meaning that $g_{|_{S^{4n-1}}}=g_{\alpha,\beta}$. By definition of the induced picture, $g$ must satisfy for all non-zero complex numbers $c$:
$$g(c\,\cdot)=\left( \frac{c}{|c|}\right)^{-\delta} \, \vert c \vert^{-i\lambda-2n} \, g.$$
\noindent We define
$$a(s,u)=\sqrt{1+4\vert s \vert ^2 + \Vert u \Vert^2}$$ 
$$r(s,u,v)= \sqrt{ a^2(s,u) + \Vert v \Vert^2}.$$
\noindent By restricting $g$ to the complex hyperplane ${\rm H} \simeq \C^{2m+1}$, we get a function $G_{\alpha,\beta}$ defined on $\C^{2m+1}$ by:
\begin{align*}
& G_{\alpha,\beta}(s,u,v) = g(1,u,2s,v)\\
& = \left( \frac{1}{r(s,u,v)} \right)^{i\lambda+2n} g_{\alpha,\beta} \left( \frac{1}{r(s,u,v)},\frac{u}{r(s,u,v)},\frac{2s}{r(s,u,v)},\frac{v}{r(s,u,v)} \right).
\end{align*}
\noindent Finally, due to the total homogeneity degree $k$ of $P_0^{\alpha,\beta}$:
$$G_{\alpha,\beta}(s,u,v)=\frac{(2s)^{\alpha}}{ \left( a^2(s,u) + \Vert v \Vert^2 \right)^{\frac{i\lambda+k}{2}+n}} \, .$$
We call $G_{\alpha,\beta}$ the \textit{non-compact form} of $P_0^{\alpha,\beta}$. The aim of the rest of Section \ref{section-non_standard_picture_and_mod._bessel_functions} is to determine the \textit{non-standard form} of $P_0^{\alpha,\beta}$, that is, $\mathcal{F}(G_{\alpha,\beta})$. We apply $\mathcal{F}_{\xi}$ in Section \ref{subsubsection-first_transform} and $\mathcal{F}_{\tau}$ in Section \ref{subsubsection-second_transform}. Calculations involve Bessel functions and various formulas that we have gathered in the appendix. 
\subsubsection{First transform} \label{subsubsection-first_transform}
\
\vskip 8pt
\noindent 
By definition (integrability condition is fulfilled):
\begin{equation} \label{equation-step_1-definition}
\displaystyle \mathcal{F}_{\xi}(G_{\alpha,\beta})(s,u,v)=\int_{\C^m}  \frac{(2s)^{\alpha}}{ \left( a^2(s,u) + \Vert \xi \Vert^2 \right)^{\frac{i\lambda+k}{2}+n}} \,e^{-2i\pi {\rm Re} \langle v,\xi \rangle}\,d\xi.
\end{equation}
\noindent In real coordinates, writing $\xi=x+iy$ and $v=a+ib$ (elements $x,y,a,b$ each belong to $\R^m$) and identifying $\xi$ and $v$ with the elements $(x,y)$ and $(a,b)$ of $\R^m \times \R^m$, formula (\ref{equation-step_1-definition}) reads:
\begin{multline}  \label{equation-step_2-real_coordinates}
\displaystyle \mathcal{F}_{\xi}(G_{\alpha,\beta})(s,u,v)=\\
\int_{\R^m \times \R^m}  \frac{(2s)^{\alpha}}{ \left( a^2(s,u) + \Vert x \Vert^2 + \Vert y \Vert^2 \right)^{\frac{i\lambda+k}{2}+n}} \,e^{-2i\pi (a \cdot x-b \cdot y)}\,dxdy.
\end{multline}
\noindent Switch to polar coordinates by writing $(x,y)=rM$ and $(a,-b)=r'M'$, with $M$ and $M'$ in $S^{2m-1} \subset \R^{2m}$, $r =\sqrt{\Vert x \Vert^2+\Vert y \Vert^2} =\Vert \xi \Vert$ and $r'=\sqrt{a^2+(-b)^2}=\Vert v \Vert$. Then the integral in Formula (\ref{equation-step_2-real_coordinates}) becomes
\begin{equation} \label{equation-step_3-polar_coordinates}
\int_0^{\infty} \left( \int_{S^{2m-1}} \frac{(2s)^{\alpha}}{ \left( a^2(s,u) + r^2 \right)^{\frac{i\lambda+k}{2}+n}} e^{-2i\pi rr'M \cdot M'} d\sigma(M) \right) r^{2m-1} dr
\end{equation}
where $M \cdot M'$ now denotes the Euclidean scalar product of $\R^{2m}$ applied to the points $M$ and $M'$ of the sphere $S^{2m-1}$ seen as vectors of $\R^{2m}$. Integral (\ref{equation-step_3-polar_coordinates}) can be written:
\begin{equation} \label{equation-step_4-polar_coordinates-2nd_version}
\displaystyle 
\int_{0}^{\infty} \frac{(2s)^{\alpha}}{ \left( a^2(s,u) + r^2 \right)^{\frac{i\lambda+k}{2}+n}} \left( \int_{S^{2m-1}} e^{-2i\pi rr'M \cdot M'} d\sigma(M) \right) r^{2m-1} dr.
\end{equation}
\noindent Proposition \ref{proposition-bochner} then changes (\ref{equation-step_4-polar_coordinates-2nd_version}) into:
\begin{equation} \label{equation-step_5-bessel_expression_1}
\int_{0}^{\infty} \frac{(2s)^{\alpha}}{ \left( a^2(s,u) + r^2 \right)^{\frac{i\lambda+k}{2}+n}} 2 \pi (rr')^{1-m} J_{m-1}(2 \pi rr') r^{2m-1} dr.
\end{equation}
\noindent Because $r'=\Vert v \Vert$, (\ref{equation-step_5-bessel_expression_1}) becomes:
\begin{equation} \label{equation-step_6-bessel_expression_2}
\displaystyle 
2^{\alpha+1} \pi s^{\alpha} \Vert v \Vert^{1-m}  \int_{0}^{\infty} \frac{r^m}{ \left( a^2(s,u) + r^2 \right)^{\frac{i\lambda+k}{2}+n}} J_{m-1}(2 \pi \Vert v \Vert r) dr.
\end{equation}
\noindent We now want to apply Proposition \ref{proposition-integral_formulas-erdelyi}. But it uses another notation system than ours. To understand how to switch from one system of notation to the other, let us define new variables $x,y,\mu$: 
$$x=r\ ;\ y=2 \pi \Vert v \Vert\ ;\ \mu=\frac{i\lambda+k}{2}+n-1\ ;\ \nu=m-1.$$
\noindent Then (\ref{equation-step_6-bessel_expression_2}) becomes:
$$2^{\alpha+1} \pi s^{\alpha} \Vert v \Vert^{1-m} y^{-\frac{1}{2}} \int_{0}^{\infty} \frac{x^{\nu+\frac{1}{2}}}{ \left( a^2(s,u) + r^2 \right)^{\mu+1}} J_{m-1}(xy) \sqrt{xy} dx.$$
\noindent Proposition \ref{proposition-integral_formulas-erdelyi} (first formula) now gives (as long as $\Vert v \Vert > 0$)
$$2^{\alpha+1} \pi s^{\alpha} \Vert v \Vert^{1-m} y^{-\frac{1}{2}} \frac{a^{\nu-\mu} y^{\mu+\frac{1}{2}} K_{\nu-\mu}(ay)}{2^{\mu} \Gamma(\mu+1)}$$
which, back to our own notation choices, is equal to
$$\frac{2^{\alpha+1} s^{\alpha} \pi^{\frac{i\lambda+k}{2}+n}}{\Gamma\left( \frac{i\lambda+k}{2}+n \right)} \left( \frac{\Vert v \Vert}{a(s,u)} \right)^{\frac{i\lambda+k}{2}+1} K_{-\left( \frac{i\lambda+k}{2}+1 \right)}(2\pi a(s,u) \Vert v \Vert).$$
\noindent To make formulas lighter, from now on we will write:
{\small $$\Lambda=\frac{i\lambda+k}{2}.$$}
Because $K_{-(\Lambda+1)}=K_{\Lambda+1}$, we have proved:
\begin{prop} \label{proposition-first_transform-k_bessel_expression}
Given any $\lambda \in \R$ and any $(k,\alpha,\beta) \in \N^3$ such that $\alpha+\beta=k$, consider the non-compact form $G_{\alpha,\beta}$ of the highest weight vector $P_0^{\alpha,\beta}$. Then for all $(s,u,v)$ in $\C \times \C^m \times \C^m$ such that $v \neq 0$:
$$\mathcal{F}_{\xi}(G_{\alpha,\beta})(s,u,v)=2^{\alpha+1} s^{\alpha} \frac{\pi^{\Lambda+n}}{\Gamma\left( \Lambda + n \right)} \left( \frac{\Vert v \Vert}{a(s,u)} \right)^{\Lambda+1} K_{ \Lambda+1 }(2\pi a(s,u) \Vert v \Vert).$$
\end{prop}
\subsubsection{Second transform} \label{subsubsection-second_transform}
\
\vskip 8pt
\noindent Consider any $(s,u,v) \in \C \times \C^m \times \C^m$ such that $s \neq 0$ and $v \neq 0$. We want to compute:
\begin{equation} \label{equation-second_transform_integral}
\mathcal{F}_{\tau} \Big( \mathcal{F}_{\xi}(G_{\alpha,\beta}) \Big) \, (s,u,v) = \int_{\C} \ \mathcal{F}_{\xi}(G_{\alpha,\beta})(\tau,u,v) \,e^{-2i\pi {\rm Re}(s\tau)} \ d\tau.
\end{equation}
\noindent Using propositions \ref{proposition-first_transform-k_bessel_expression} and \ref{proposition-bessel_asymptotics-erdelyi}, one shows  that $\tau \longmapsto \mathcal{F}_{\xi}(G_{\alpha,\beta})(\tau,u,v)$ is indeed integrable. Let us use letters $a,b,x,y$ again, this time taking them to refer to real numbers: $s=a+ib$ and $\tau=x+iy$. Then ${\rm Re}(s\tau)=ax-by$ and (\ref{equation-second_transform_integral}) becomes
$$\int_{\R^2} \frac{2^{\alpha+1} \tau^{\alpha} \pi^{\Lambda+n}}{\Gamma\left( \Lambda + n \right)} \left( \frac{\Vert v \Vert}{a(\tau,u)} \right)^{\Lambda+1} K_{ \Lambda+1 }(2\pi a(\tau,u) \Vert v \Vert) \,e^{-2i\pi (ax-by)}\,dxdy$$
which can be re-organised as
\begin{equation} \label{equation-ax_by_integral}
\frac{2^{\alpha+1} \pi^{\Lambda+n} \Vert v \Vert^{\Lambda+1}}{\Gamma\left( \Lambda + n \right)} \int_{\R^2} \frac{\tau^{\alpha} K_{ \Lambda+1 }(2\pi a(\tau,u) \Vert v \Vert) }{\Big( a(\tau,u) \Big)^{\Lambda+1}} \,e^{-2i\pi (ax-by)}\,dxdy.
\end{equation}
\noindent Let us again use polar coordinates (outside the origin):
\begin{itemize}
\item $(x,y)= r v_{\theta}$ with $r>0$, $\theta \in \R$ and $v_{\theta}=(\cos \theta, \sin \theta)$; thus, $\tau =r e^{i\theta}$.
\item $(a,-b)=r' v_{\theta'}$ with $r'>0$, $\theta' \in \R$ and $v_{\theta'}=(\cos \theta', \sin \theta')$; thus, $\overline{s} =r' e^{i\theta'}$.
\end{itemize}
\noindent Let us write $a(r,u)$ instead of $a(\tau,u)$:
$$a(r,u)=\sqrt{1+4r^2+\Vert u \Vert^2}.$$
\noindent Integral (\ref{equation-ax_by_integral}) can now be written:
\begin{multline} \label{equation-cos_sin_integral}
\frac{2^{\alpha+1} \pi^{\Lambda+n} \Vert v \Vert^{\Lambda+1}}{\Gamma\left( \Lambda + n \right)} \int_0^{\infty} \frac{r^{\alpha} K_{ \Lambda+1 }(2\pi a(r,u) \Vert v \Vert) }{\Big( a(r,u) \Big)^{\Lambda+1}} \\
\left( \int_0^{2\pi} e^{i\alpha\theta}\,e^{-2i\pi rr' (\cos \theta \cos \theta'+\sin \theta \sin \theta')} d\theta \right)\,rdr.
\end{multline}
\noindent Following Proposition \ref{proposition-bessel_integral_formula}, the inner integral
$$\int_0^{2\pi} e^{i\alpha \theta} \, e^{-2i\pi r r' (\cos \theta \, \cos \theta' \, + \, \sin \theta \, \sin \theta')} \, d\theta$$
is equal to:
\begin{equation} \label{equation-inner_integral_computation}
2 \pi e^{i\alpha\left( \theta'-\frac{\pi}{2} \right) } \, J_{\alpha}(2\pi r r').
\end{equation}
\noindent Because $r'=\vert s \vert$ and $\theta'={\rm Arg}(\overline{s})$, (\ref{equation-inner_integral_computation}) is equal to:
$$2 \pi e^{i\alpha \left( {\rm Arg}(\overline{s}) -  \frac{\pi}{2} \right) } \, J_{\alpha}(2\pi r \vert s \vert ).$$
\noindent This turns (\ref{equation-cos_sin_integral}) into:
\begin{multline} \label{equation-non_standard_version_of_hwv_gamma=0}
\frac{2^{\alpha+2} \pi^{\Lambda+n+1} \Vert v \Vert^{\Lambda+1} e^{i\alpha \left( {\rm Arg}(\overline{s}) -  \frac{\pi}{2} \right) } }{\Gamma\left( \Lambda + n \right)} \\
\int_0^{\infty} \frac{r^{\alpha + 1} \, K_{ \Lambda+1 }(2\pi a(r,u) \Vert v \Vert) }{\Big( a(r,u) \Big)^{\Lambda+1}} \, J_{\alpha}(2\pi r |s|) \, dr.
\end{multline}
\noindent We can now apply the second formula of proposition \ref{proposition-integral_formulas-erdelyi}. To help follow notation choices made in this proposition, we set:
\begin{itemize}
\item $x=2r$, $dx=2dr$ and $\beta=\sqrt{1+\Vert u \Vert^2 } > 0$;
\item $a=2\pi\Vert v \Vert > 0$ (careful: this variable $a$ is not what we have denoted $a(r,u)$);
\item $y=\pi\vert s\vert > 0$, $\nu=\alpha$ and $\mu=\Lambda+1$.
\end{itemize}
\noindent Plugging these expressions in (\ref{equation-non_standard_version_of_hwv_gamma=0}) and using Proposition \ref{proposition-integral_formulas-erdelyi}, we finally obtain:
\begin{thm} \label{theorem-final_formula}
Given any $\lambda \in \R$ and any $(k,\alpha,\beta) \in \N^3$ such that $\alpha+\beta=k$, consider the non-compact form $G_{\alpha,\beta}$ of the highest weight vector $P_0^{\alpha,\beta}$. Then for all $(s,u,v)$ in $\C \times \C^m \times \C^m$ such that $s \neq 0$ and $v \neq 0$:
\begin{multline*}
\mathcal{F}(G_{\alpha,\beta})(s,u,v)\, = \\
\int_{\C \times \C^m} \, \frac{(2\tau)^{\alpha}}{ \left(1+4\vert \tau \vert ^2 + \Vert u \Vert^2 + \Vert \xi \Vert^2 \right)^{\frac{i\lambda+k}{2}+n}} \, e^{-2i\pi  {\rm Re} \left( s \tau + \langle v,\xi \rangle \right) } \,d\tau d\xi \, = \\
R(s,u,v) \, K_{\frac{i\lambda+\delta}{2}} \left( \pi \sqrt{ 1 + \Vert u \Vert^2 } \sqrt{ \vert s \vert^2+4\Vert v \Vert^2 } \right)
\end{multline*}
where we set
$$R(s,u,v) \, = \, \frac{(-i \, \overline{s})^{\alpha} \, \pi^{i\lambda+\beta+n}}{2^{\frac{i\lambda+k}{2}+1} \, \Gamma \left( \frac{i\lambda+k}{2}+n \right)} \left( \frac{ \sqrt{ |s|^2+4\Vert v \Vert^2}}{\pi \sqrt{1+\Vert u \Vert^2}} \right)^{\frac{i\lambda+\delta}{2}}.$$
\end{thm}
\noindent This theorem establishes the second part of our "Main results" theorem (stated in the introduction).
\begin{nnsrem} \
\begin{itemize}
\item The partial Fourier transform $\mathcal{F}$ changes, up to a constant $\frac{-1}{i\pi}$, multiplication by the $s$ coordinate into differentiation with respect to $s$. This implies $\mathcal{F}(G_{\alpha+1,\beta-1})=\frac{2}{-i\pi}\frac{\partial}{\partial s} \Big( \mathcal{F}(G_{\alpha,\beta}) \Big)$, which is in fact proved in \cite{Mendousse2017} for $\beta \geq 2$ (due to integrability issues).
\item One inevitably notices in the formula of Theorem \ref{theorem-final_formula} a structure in two variables. Indeed, the particular value $\alpha=0$ and the square root terms lead us to study the functions
$$\begin{array} {cccc}
\psi_{\nu} : & ]0,+\infty[ \times ]0,+\infty[ & \longrightarrow & \C \\
         & (x,y)           & \longmapsto & \left( \frac{x}{y} \right)^{\nu} K_{\nu} (xy) \\
\end{array}$$
where the parameter $\nu$ is a complex number. Fix $\nu \in \C$, $x_0 > 0$ and define the function
$$\begin{array} {cccc}
\varphi_{x_0,\nu} : & ]0,+\infty[ & \longrightarrow & \C \\
         & y           & \longmapsto & \psi_{\nu}(x_0,y). \\
\end{array}$$
\noindent One can show:
$$ \varphi_{x_0,\nu}''(y) \, + \, \frac{(1+2\nu)}{y} \, \varphi_{x_0,\nu}'(y) \, - \, x_0^2 \, \varphi_{x_0,\nu}(y) \, = \, 0.$$
\noindent This equation belongs to the family of  \textit{Emden-Fowler equations} (or \textit{Lane-Emden equations}) and its solutions can be written as the following combinations of Bessel functions of the first and second kind:
$$u(t) \, = \, C_1 \, t^{-\nu} \, J_{\nu} \left( -itx_0 \right) \, + \, C_2 \, t^{-\nu} \, Y_{\nu} \left( -itx_0 \right).$$
Similar conclusions hold if one fixes $y_0$ instead of $x_0$.
\end{itemize}
\end{nnsrem}
\section{Appendix: Bessel functions} \label{section-appendix}
\noindent In these definitions, following for instance \cite{Lebedev1972} (sections 5.3 and 5.7), we take $\nu \in \C$ and $z \in \C \, \setminus \{0\}$ such that $-\pi < {\rm Arg}(z) < \pi$:
\begin{enumerate}
\item The \textit{Bessel function of the first kind} is the function $J_{\nu}$ defined by:
$$J_{\nu}(z)=\sum_{k=0}^{\infty} \frac{(-1)^k}{\Gamma(k+1) \Gamma(k+\nu+1)} \left( \frac{z}{2} \right)^{\nu+2k}.$$
\item The \textit{Bessel function of the second kind} is the function $Y_{\nu}$ defined by:
$$Y_{\nu}(z)=\frac{J_{\nu}(z) \cos(\nu \pi) - J_{-\nu}(z)}{\sin (\nu \pi)}$$
when $\nu \notin \Z  $ and, when $\nu \in \Z$, by
$$Y_{\nu}(z)=\lim_{\substack{
\epsilon \rightarrow \nu \\
0< |\epsilon-\nu| < 1}
} Y_{\epsilon}(z).$$
\item The \textit{modified Bessel function of the first kind} is the function $I_{\nu}$ defined by: 
$$I_{\nu}(z)=\sum_{k=0}^{\infty} \frac{1}{\Gamma(k+1) \Gamma(k+\nu+1)} \left( \frac{z}{2} \right)^{\nu+2k}.$$
\item The \textit{modified Bessel function of the third kind} is the function $K_{\nu}$ defined by 
$$K_{\nu}(z)=\frac{\pi}{2} \frac{I_{-\nu}(z) - I_{\nu}(z)}{\sin (\nu \pi)}$$
when $\nu \notin \Z  $ and, when $\nu \in \Z$, by
$$K_{\nu}(z)=\lim_{\substack{
\epsilon \rightarrow \nu \\
0 < |\epsilon-\nu| < 1}
} K_{\epsilon}(z).$$
\end{enumerate}
\noindent As a consequence of Formula (2) of Section 7.3.1 in Chapter VII of \cite{Erdelyi-higher_etc.-2}:
\begin{prop}[An integral representation of Bessel functions] \label{proposition-bessel_integral_formula} Given $\nu \in \N$, $\rho >0$ and $a > 0$:
$$J_{\nu}(\rho) \, = \, \frac{1}{2\pi e^{i\nu\left( a-\frac{\pi}{2} \right)}} \, \int_0^{2\pi} e^{i\nu \theta} \, e^{-i\rho(\cos a \, \cos \theta \, + \, \sin a \, \sin \theta)} \, d\theta.$$
\end{prop}
\noindent Formulas in the next proposition are stated in \cite{Erdelyi-tables_etc.-2} (Chapter VIII: Formula (20) of Section 8.5 and Formula (35) of Section 8.14):
\begin{prop}[Two integral formulas involving Bessel functions] \label{proposition-integral_formulas-erdelyi}\
\begin{itemize}
\item For any real number $y > 0$ and any complex numbers $a,\nu,\mu$ such that ${\rm Re}(a) > 0$ and $-1 < {\rm Re}(\nu) < 2 {\rm Re}(\mu) + \frac{3}{2}$, one has:
$$\int_0^{\infty} x^{\nu + \frac{1}{2}} \left( x^2 + a^2 \right)^{-\mu-1} J_{\nu}(xy) \sqrt{xy} \ dx =$$
$$ \frac{a^{\nu-\mu} y^{\mu+\frac{1}{2}} K_{\nu-\mu}(ay)}{2^{\mu} \Gamma(\mu+1)}.$$
\item For any real number $y > 0$ and any complex numbers $a,\beta,\nu,\mu$ such that ${\rm Re}(a) > 0$, ${\rm Re}(\beta) > 0$ and ${\rm Re}(\nu) > -1$, one has:
$$\int_0^{\infty} x^{\nu + \frac{1}{2}} \left( x^2 + \beta^2 \right)^{-\frac{\mu}{2}} K_{\mu} \left( a(x^2+\beta^2)^{\frac{1}{2}} \right)  J_{\nu}(xy) \sqrt{xy}\ dx =$$
$$a^{-\mu} \beta^{\nu+1-\mu} y^{\nu + \frac{1}{2}} (a^2+y^2)^{\frac{\mu}{2}-\frac{\nu}{2}-\frac{1}{2}} K_{\mu-\nu-1} \left( \beta (a^2+y^2)^{\frac{1}{2}} \right).$$
\end{itemize}
\end{prop}
\noindent One can find the next formula in \cite{Erdelyi-higher_etc.-2} (Section 7.4.1, Formula (4)):
\begin{prop}[Asymptotic expansion for modified Bessel functions] \label{proposition-bessel_asymptotics-erdelyi} 
For any fixed $P \in \N \setminus \{0\}$ and $\nu \in \C$:
$$K_{\nu}(z)=\left( \frac{\pi}{2z} \right)^{\frac{1}{2}} e^{-z} \left( \left[ \sum_{p=0}^{P-1} \frac{ \Gamma \left( \frac{1}{2} + \nu + p\right) }{p! \ \Gamma \left( \frac{1}{2} + \nu - p\right) } \ (2z)^{-p} \right] + {\rm O} \left( |z|^{-P} \right) \right).$$
\end{prop}
\noindent The following proposition can be derived from Section 9.6 in \cite{Faraut2008}:
\begin{prop}[Bochner formula] \label{proposition-bochner} 
Consider any integer $p \geq 2$. For $\xi'\in S^{p-1}$ and $s > 0$:
$$\int_{S^{p-1}}e^{-2i\pi s \xi \cdot \xi'} d\sigma(\xi)
=2 \pi s^{1-\frac{p}{2}} J_{\frac{p}{2}-1}(2\pi s)$$
where $d\sigma$ denotes the Euclidean measure of $S^{p-1}$ and $\xi \cdot \xi'$ denotes the Euclidean scalar product of $\R^p$ applied to the elements of the sphere $\xi$ and $\xi'$ seen as elements of $\R^p$.
\end{prop}

\footnotesize{ 
\noindent \textbf{Contact information.}\\
\noindent Universit{\'e} de Reims Champagne-Ardenne\\
\noindent Laboratoire de Math{\'e}matiques FRE 2011 CNRS\\
\noindent UFR Sciences Exactes et Naturelles, Moulin de la Housse - BP 1039\\
\noindent 51687 REIMS Cedex 2, FRANCE\\
\noindent \texttt{gregory.mendousse@univ-reims.fr}
}

\end{document}